\documentclass{amsart}

\usepackage{epsfig,amsmath,amssymb}

\begin{document}

\bibliographystyle{siam}

\title{Kazhdan-Lusztig polynomials for\\ 321-hexagon-avoiding permutations}

\author{Sara C. Billey}
\email{billey@math.mit.edu} 

\address{Author's address: 
Dept.\ of Mathematics, 2-363c\\ 
Massachusetts Institute of Technology\\ Cambridge, MA 02139}

\author{Gregory S. Warrington}
\email{gwar@math.harvard.edu} 

\address{Author's address: 
Dept.\ of Mathematics\\
Harvard University\\ Cambridge, MA 02138}

\subjclass{05E15 (Primary); 20F55, 32S45, 14M15 (Secondary)}

\begin{abstract}
  In \cite{Deodhar90}, Deodhar proposes a combinatorial framework for
  determining the Kazhdan-Lusztig polynomials $P_{x,w}$ in the case
  where $W$ is any Coxeter group.  We explicitly describe the
  combinatorics in the case where $W=\mathfrak{S}_n$ (the symmetric
  group on $n$ letters) and the permutation $w$ is
  321-hexagon-avoiding.  Our formula can be expressed in terms of a
  simple statistic on all subexpressions of any fixed reduced
  expression for $w$.  As a consequence of our results on
  Kazhdan-Lusztig polynomials, we show that the Poincar\'e polynomial
  of the intersection cohomology of the Schubert variety corresponding
  to $w$ is $(1+q)^{l(w)}$ if and only if $w$ is 321-hexagon-avoiding.
  We also give a sufficient condition for the Schubert variety $X_w$
  to have a small resolution.  We conclude with a simple method for
  completely determining the singular locus of $X_w$ when $w$ is
  321-hexagon-avoiding.  The results extend easily to those Weyl
  groups whose Coxeter graphs have no branch points ($B_n$, $F_4$,
  $G_2$).
\end{abstract}

\maketitle

\begin{center}
  \textit{In memory of Rodica Simion}
\end{center}


\renewcommand{\sectionmark}[1]{\markright{\thesection. #1}}



\newtheorem{thm}{Theorem}
\newtheorem{lem}{Lemma}
\newtheorem{cor}{Corollary}
\newtheorem{prop}{Proposition}
\newtheorem{exmp}{Example}
\newtheorem{REM}{Remark}


\newtheorem{dfn}{Definition}
\newtheorem{fact}{Fact}

\newcommand{\thmref}[1]{Theorem~\ref{#1}}
\newcommand{\corref}[1]{Corollary~\ref{#1}}
\newcommand{\propref}[1]{Proposition~\ref{#1}}
\newcommand{\lemref}[1]{Lemma~\ref{#1}}
\newcommand{\lemsref}[1]{Lemmas~\ref{#1}}
\newcommand{\factref}[1]{Fact~\ref{#1}}
\newcommand{\rmkref}[1]{Remark~\ref{#1}}
\newcommand{\defref}[1]{Definition~\ref{#1}}
\newcommand{\egref}[1]{Example~\ref{#1}}
\newcommand{\figref}[1]{Figure~\ref{#1}}
\newcommand{\secref}[1]{Section~\ref{#1}}

\newcommand{\sfrac}[2]{\genfrac{\{}{\}}{0pt}{}{#1}{#2}}
\newcommand{\sumsb}[1]{\sum_{\substack{#1}}}  
        
\newcommand{\br}{\mathbf{a}}
\newcommand{\brp}{\mathbf{a'}}
\newcommand{\bsig}{{\boldsymbol{\sigma}}}
\newcommand{\bnu}{{\boldsymbol{\nu}}}
\newcommand{\bgam}{{\boldsymbol{\gamma}}}
\newcommand{\bdel}{{\boldsymbol{\delta}}}

\newcommand{\rank}{\operatorname{lvl}}
\newcommand{\pt}{\operatorname{pt}}
\newcommand{\heap}{\operatorname{Heap}}
\newcommand{\lcz}{\operatorname{lcz}}
\newcommand{\rcz}{\operatorname{rcz}}
\newcommand{\mcz}{\operatorname{mcz}}
\newcommand{\ver}{\operatorname{ver}}
\newcommand{\codim}{\operatorname{codim}}
\newcommand{\bcone}{\operatorname{Cone_{\wedge}}}
\newcommand{\ucone}{\operatorname{Cone^{\vee}}}

\newcommand{\IH}{\operatorname{IH}}

\newcommand{\defeq}{\overset{\text{def}}{=}}

\newcommand{\io}{{i_1}}
\newcommand{\iw}{{i_2}}
\newcommand{\inn}{{i_r}}

\newcommand{\bbC}{\mathbb{C}}
\newcommand{\bbZ}{\mathbb{Z}}
\newcommand{\bbN}{\mathbb{N}}
\newcommand{\bbQ}{\mathbb{Q}}

\newcommand{\cP}{{\mathcal P}}
\newcommand{\xsing}{{X_w^{\text{sing}}}}

\newcommand{\cpw}{{\mathcal P}(\br)}
\newcommand{\cpww}{{\mathcal P}_w(\br)}
\newcommand{\cpws}{{\mathcal P}(\br s)}
\newcommand{\cpxw}{{\mathcal P}_x(\br)}
\newcommand{\cpyw}{{\mathcal P}_y(\br)}
\newcommand{\cpxsw}{{\mathcal P}_{xs}(\br)}
\newcommand{\cpxsws}{{\mathcal P}_{xs}(\br/s)}
\newcommand{\cpxws}{{\mathcal P}_x(\br/s)}
\newcommand{\cpzero}{{\mathcal P}_x^0(\br)}
\newcommand{\cpone}{{\mathcal P}_x^1(\br)}
\newcommand{\cpeps}{{\mathcal P}_x^\epsilon(\br)}

\newcommand{\cd}{{\mathcal D}}
\newcommand{\Dbr}{\Delta_\bsig}
\newcommand{\Dbrj}{\Delta_{\bsig[j]}}
\newcommand{\Dbrr}{\Delta_{\bsig[r]}}
\newcommand{\Dbrk}{\Delta_{\bsig[k]}}
\newcommand{\Dbrjm}{\Delta_{\bsig[j-1]}}

\newcommand{\pw}{P(\br)}
\newcommand{\pww}{P_w(\br)}
\newcommand{\pxw}{P_x(\br)}
\newcommand{\pxsw}{P_{xs}(\br)}
\newcommand{\pxws}{P_x(\br/s)}
\newcommand{\pxsws}{P_{xs}(\br/s)}

\newcommand{\pzw}{P_z(\br)}
\newcommand{\pzsw}{P_{zs}(\br)}
\newcommand{\pzws}{P_z(\br/s)}
\newcommand{\pzsws}{P_{zs}(\br/s)}

\newcommand{\pew}{P_e(\br)}
\newcommand{\cpew}{{\mathcal P}_e(\br)}

\newcommand{\rxz}{R_{x,z}}
\newcommand{\rxsz}{R_{xs,z}}
\newcommand{\rxzs}{R_{x,zs}}
\newcommand{\rxszs}{R_{xs,zs}}

\newcommand{\sij}{s_{i_j}}
\newcommand{\sik}{s_{i_k}}


\section{Introduction}\label{sec:intro}

In \cite{K-L1}, Kazhdan and Lusztig constructed certain
representations of the Hecke algebra associated to a Coxeter group $W$
in order to elucidate representation-theoretic questions concerning
$W$ itself.  To do this, they introduced a class of polynomials now
known as the Kazhdan-Lusztig polynomials.  These polynomials were
quickly seen to play an important role in Lie theory.  For instance,
they give a natural setting for expressing multiplicities of
Jordan-H\"older series of Verma modules (see \cite{BeilBern,BryKash}).
Introductions to these polynomials can be found in
\cite{Brenti98,Deodhar94,Hum}.

While there are many interpretations of, and uses for, these
polynomials, their combinatorial structure is far from clear.  Kazhdan
and Lusztig originally defined the polynomials in terms of a
complicated recursion relation.  In \cite{K-L1}, it was conjectured
that the coefficients of these polynomials are non-negative.  This has
been proved for many important $W$ (such as (affine) Weyl groups)
\cite{K-L2}, but not for arbitrary Coxeter groups.  There has been
limited success in finding non-recursive formulas for the
Kazhdan-Lusztig polynomials.  Brenti \cite{Brenti94,Brenti97} has
given a non-recursive formula in terms of an alternating sum over
paths in the Bruhat graph.  Lascoux and Sch\"utzenberger \cite{LS11}
have given an explicit formula for $P_{x,w}$ in the case where $W$ is
the symmetric group and $x,w$ are Grassmannian permutations.
Zelevinsky \cite{zel} has even constructed a small resolution of $X_w$
in this case.  Lascoux \cite{Lascoux95} extends the results of
\cite{LS11} to twisted vexillary permutations.  Finally, Shapiro,
Shapiro and Vainshtein \cite{ssv} and Brenti and Simion \cite{Bren-Sim}
find explicit formulas for certain classes of permutations.

Deodhar \cite{Deodhar90} proposes a combinatorial framework for
determining the Kazhdan-Lusztig polynomials for an arbitrary Coxeter
group.  The algorithm he describes is shown to work for all Weyl
groups.  However, the algorithm is impractical for routine
computations.  In this paper, we utilize Deodhar's framework to
calculate $P_{x,w}$ for 321-hexagon-avoiding elements $w\in
\mathfrak{S}_n$.  For these elements, Deodhar's algorithm turns out to be
trivial.  As a result, in these cases we get a very explicit
description of the polynomials.  The algorithm consists of calculating
Deodhar's defect statistic on each subexpression of a given reduced
expression.  We also show that the property of $w$ being
321-hexagon-avoiding is equivalent to several nice properties on $w$
in the Hecke algebra and in the cohomology of the corresponding
Schubert variety $X_w$.  In particular, we have the following (the
necessary definitions can be found in \secref{sec:prelim} and
\secref{sec:hex}):

\begin{thm}\label{mainthm}
  Let $\br = s_{i_1}\cdots s_{i_r}$ be a reduced expression for $w\in
  \mathfrak{S}_n$.  The following are equivalent:
  \begin{enumerate}
  \item $w$ is 321-hexagon-avoiding.\label{itone}
  \item Let $P_{x,w}$ denote the Kazhdan-Lusztig polynomial for $x \leq
  w$.  Then 
    \begin{equation}
      P_{x,w}= \sum q^{d(\sigma)}
    \end{equation}
    where $d(\sigma)$ is the defect statistic and the sum is over all
    masks $\sigma$ on $\mathbf{a}$ whose product is $x$.\label{ittwo}
  \item \label{itthree} The Poincar\'e polynomial for the full intersection 
    cohomology   group of $X_w$ is 
    \begin{equation}
      \sum_i \dim(\IH^{2i}(X_w))q^i = (1+q)^{l(w)}.
    \end{equation}
  \item The Kazhdan-Lusztig basis element $C_w'$ satisfies $C_w' =
    C_{s_{i_1}}'\cdots C_{s_{i_r}}'$.\label{itfour}
  \item The Bott-Samelson resolution of $X_w$ is small.\label{itfive}
  \item $\IH_*(X_w) \cong H_*(Y)$, where $Y$ is the Bott-Samelson
    resolution of $X_w$.\label{itsix}
  \end{enumerate}
\end{thm}

\begin{REM}
  Equivalence of \ref{ittwo}, \ref{itfour} and \ref{itfive} is
  implicit in Deodhar \cite{Deodhar90}.
\end{REM}

\begin{REM}
  Lusztig \cite{lusztig} and Fan and Green \cite{FanGreen} have
  already studied those elements $w$ for which part \ref{itfour} of
  the main theorem hold.  In the terminology of these papers, such a
  $w$ is ``tight.''  Also, Fan and Green show the implication
  \textit{\ref{itfour} $\Longrightarrow$ \ref{itone}} of \thmref{mainthm}.
\end{REM}

\begin{REM}
  For concreteness, this paper refers only to $\mathfrak{S}_n$.
  However, 2 through 6 hold for all Weyl groups.  In addition, our
  combinatorial characterization of \ref{itone} $\Longleftrightarrow$
  \ref{ittwo} can be extended to the other ``non-branching'' Weyl
  groups $B_n,F_4,G_2$ (see \cite{gwar}).  One need simply replace
  ``321-avoiding'' by ``short-braid-avoiding'' in any statements made
  (e.g., ``321-hexagon-avoiding'' $\mapsto$
  ``short-braid-hexagon-avoiding'').  The characterization in
  \ref{itone} fails to hold for $D_n,E_6,E_7,E_8$ primarily due to
  failure of \lemref{latcon}.  An appropriate analogue of
  hexagon-avoiding for these other Weyl groups would fix this
  deficiency.
\end{REM}

The organization of the paper is as follows.  In \secref{sec:prelim}
we introduce necessary background definitions.  In \secref{sec:hex} we
introduce the notion of pattern avoidance and in \secref{sec:frame} we
present Deodhar's combinatorial framework.  A critical tool used to
prove \thmref{mainthm} is the defect graph explored in
\secref{sec:graph}.  In \secref{sec:pf} this graph is used to prove
\thmref{mainthm}.  \secref{sec:conj} contains an application of
\thmref{mainthm} to a conjecture of Haiman.  \secref{sec:loci}
determines the singular locus of Schubert varieties corresponding to
321-hexagon-avoiding permutations.  Finally, \secref{sec:enum}
contains a table enumerating the elements of $\mathfrak{S}_n$ for which
\thmref{mainthm} applies.  We do not know a closed form for this
sequence.

\section{Preliminaries}\label{sec:prelim}

Let $\mathfrak{S}_n$ denote the symmetric group on $n$ letters.
Choose the standard presentation $\mathfrak{S}_n = \langle
s_1,\ldots,s_{n-1} : s_i^2 = 1, s_is_j = s_j s_i$ for $|i-j| > 1$, and
$s_is_{i+1} s_i = s_{i+ 1} s_i s_{i+1} \rangle.$ Let $S =
\{s_i\}_{i\in [1\ldots n-1]}$ denote the generating set for
$\mathfrak{S}_n$.  An \textit{expression} is any product of generators
$s_i$.  The \textit{length} $l(w)$ of an element $w\in \mathfrak{S}_n$
is the minimum $r$ for which we have an expression $w = s_\io\cdots
s_{i_r}$.  A \textit{reduced expression} $w=s_{i_1}\cdots s_{i_r}$ is
an expression for which $l(w) = r$.  If $v,w\in \mathfrak{S}_n$, then
$v \leq w$ will signify that $v$ is below $w$ in the Bruhat-Chevalley
order (see, e.g., \cite{Hum}).  This order is characterized by $v\leq
w$ if and only if every reduced expression for $w$ contains a
subexpression for $v$.

For the remainder of this section, all of our definitions apply to any
finite Weyl group $W$.  However, following this section, we will
restrict our attention to the case where $W = \mathfrak{S}_n$.

In order to define the Kazhdan-Lusztig polynomials, we now recall the
notion of the Hecke algebra $\mathcal{H}$ associated to a finite Weyl
group $W$.  $\mathcal{H}$ has basis $T_w$ indexed by the elements of
$W$.  For all generators $s$ of $W$, we have
\begin{align}
  \label{defrel}
  T_s T_w &= T_{s w} \text{ if } l(sw) > l(w), \\
  T_s^2 &= (q-1)T_s + q T_e
\end{align}
(where $e$ is the identity element of $W$).  This is an algebra over
$A = \bbQ(q^{1/2})$.  Following \cite{K-L1}, we define an involution
on $A$ by $\overline{q^{1/2}} = q^{-1/2}$.  Extend this to an
involution on $\mathcal{H}$ by setting
\begin{equation}
  \iota(\sum_w \alpha_w T_w) =
  \sum_w \bar{\alpha}_w (T_{w^{-1}})^{-1}.
\end{equation}
From \cite{K-L1}, we have that the Kazhdan-Lusztig polynomials are
determined uniquely by the following:

\begin{thm}[Theorem 1.1, \cite{K-L1}]\label{klthm}
  For any $w\in W$, there is a unique element $C_w'\in \mathcal{H}$
  such that 
  \begin{enumerate}
  \item \ $C_w' = q^{-l(w)/2} \sum_{x\leq w} P_{x,w} T_x$, and
  \item \ $\iota(C_w') = C_w'$,
  \end{enumerate}
  where $P_{x,w}\in A$ is a polynomial in $q$ of degree at most
  $\frac{1}{2}(l(w)-l(x)-1)$ for $x < w$, $P_{w,w} = 1$, and $P_{x,w}
  = 0$ if $x \not\leq w$.
\end{thm}
  
As mentioned above, it is conjectured in \cite{K-L1} that the
coefficients of $P_{x,w}$ are non-negative.

Several of the conditions in \thmref{mainthm} require some notation
regarding cohomology.  So let $W$ be the Weyl group of some semi-simple
algebraic group $G$ with Borel subgroup $B$.  $C_w$ will denote the
Schubert cell in the flag variety $G/B$ corresponding to $w\in W$
(see, e.g., \cite{bou}).  $X_w$ will denote the corresponding Schubert
variety, $X_w = \cup_{v\leq w} C_v$.  For any variety $X$ (such as
some $X_w$), we let $\IH^i(X)$ denote the $i$-th (middle) intersection
cohomology group of $X$.  Suppose that $f: Y \longrightarrow X$ is a
resolution of singularities of $X$.  The map $f$ is said to be a
\textit{small resolution} if for every $r > 0$,
\begin{equation}
  \codim\{x\in X: \dim f^{-1}(x) \geq r \} > 2r.
\end{equation}
A commonly used resolution of the singularities of $X_w$ is the
Bott-Samelson resolution (see \cite{B-S,Dem}).  \thmref{mainthm}
yields an easy criterion for determining when such a resolution is
small.

\section{Pattern Avoidance and Heaps}\label{sec:hex}

It will be useful to view elements of $\mathfrak{S}_n$ as permutations on
$[1,2,\ldots,n]$.  To this end, we identify $s_i$ with the
transposition $(i,i+1)$.  Let $w(i)$ be the image of $i$ under the
permutation $w$.  Hence, we have a one-line notation for a permutation
$w$ given by writing the image of $[1,2,\ldots,n]$ under the action of
$w$: $[w(1),w(2),\ldots,w(n)]$.

The results of this paper pertain to a particular set of elements of
$\mathfrak{S}_n$.  This subset will be defined using the notion of pattern
avoidance.  Let $v\in S_k$ and $w\in S_l$.  Say that $w$
\textit{avoids} $v$ (or is $v$-avoiding) if there do not exist $1\leq
i_1<\cdots < i_k\leq l$ with $w(i_1),w(i_2),\ldots,w(i_k)$ in the same
relative order as $v(1),v(2),\ldots,v(k)$.  We are interested in two
particular instances of pattern avoidance.  The first is where $v =
[3,2,1]$.  It is shown by Billey-Jockusch-Stanley \cite{BJS} that the
321-avoiding permutations in $\mathfrak{S}_n$ are precisely those for which no
reduced expression contains a substring of the form $s_is_{i\pm
1}s_i$.  In the context of reduced expressions, 321-avoiding
permutations are called \textit{short-braid-avoiding} (terminology due
to Zelevinsky, according to \cite{Fan1}).  Short-braid-avoiding
permutations have been studied by Fan and Stembridge
\cite{fan,Fan1,Stem4,Stem3}.

The second instance of pattern-avoidance with which we will be
concerned is most easily visualized via a poset associated to $w$.  So
let $w\in \mathfrak{S}_n$ be 321-avoiding and fix some reduced expression
$\br = s_{i_1}\cdots s_{i_r}$ for $w$.  By \cite{T}, all reduced
expressions for such a 321-avoiding $w$ are equivalent up to moves of
the form $s_is_j \rightarrow s_js_i$ for $|i-j| > 1$.  This allows us
to associate a well-defined poset to $w$ (rather than just to $\br$,
see \cite{Stem4}).  Let the generators $\{\sij\}_{j=1}^r$ in our
reduced expression label the elements of our poset.  For an ordering,
we take the transitive closure of
\begin{quote}
  $\sij \preccurlyeq \sik$ if $s_{i_{j+1}}\dots s_{i_{k-1}}\sik= \sik s_{i_{j+1}}\dots s_{i_{k-1}}$ and  $\sij\sik \neq \sik \sij$.
\end{quote}
We now wish to embed this poset in the plane in a very particular way.
Effectively, what we do is send a generator $\sij$ to the point in the
plane $(i_j,\rank(\sij))$ where $\rank(\sij)$ measures the maximal length
of a chain $s_{i_b} \preccurlyeq \ldots \preccurlyeq \sij$ over all $b
\leq j$.  However, in order for our embedding to have the properties
we need, this procedure needs to be adjusted slightly.

So, as above, embed this poset in the plane via $\sij \mapsto \pt(j)
\defeq (i_j,\rank(\sij))$, where we define $\rank(\sij)$ as
follows: Let $k$ be as small as possible in the interval
$[1,\ldots,j]$ such that $s_{i_j}$ commutes with $s_{i_l}$ for all $l$
with $k\leq l \leq j$.  Now, initially, define a level function by:
$\rank_L(s_{i_j}) = 0$ if $k=1$ and $\rank_L(\sij) =
\rank_L(s_{i_{k-1}}) + 1$ if $k\geq 2$.

\begin{figure}[htbp]
  \begin{center}
    \leavevmode
    \input{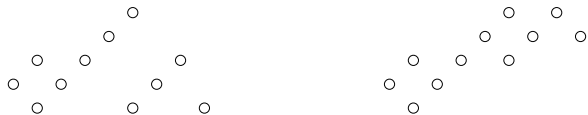}
\caption{Let $w = s_9s_6s_7s_8s_2s_1s_3s_2s_4s_5s_6$.  The left image shows the result
  of the embedding $\sij \mapsto (i_j,\rank_L(\sij))$.  On the right
  is the result of pushing the ``connected components'' together.}
    \label{fig:poset}
  \end{center}
\end{figure}

For most purposes, $\rank_L(\cdot)$ gives us what we'd like. However,
with $\rank_L(\cdot)$ as the level function, ``connected components''
do not necessarily abut.  \figref{fig:poset} gives an example of the
embedding $(i_j,\rank_L(\sij))$ and how it can be improved by
coalescing ``connected components.''

So, we first define connected components by imposing an equivalence
$\sim$ on the generators in our expression for w: Let $s_{i_j} \sim
s_{i_k}$ if $i_j = i_k\pm 1$ and $\rank_L(s_{i_j}) =
\rank_L(s_{i_k})\pm 1$.  Extend this equivalence transitively.  Now,
since we are assuming that $w$ is 321-avoiding, the components have a
canonical partial order.  It is then a simple matter to uniformly
adjust the levels of all members of a particular connected component to
allow distinct components to abut as much as possible and hence
``coalesce.''  Define $\rank(\sij)$ to be this adjustment of the level
$\rank_L(\sij)$.

We will refer to the realization $\sij \mapsto (i_j,\rank(\sij)$ of
our poset as \textit{$\heap(w)$}.  The notion of $\heap(w)$ is due to
Viennot \cite{Viennot}, see also the work of Stembridge \cite{Stem4}
in the context of fully-commutative elements.  Note that $\sij$ can
cover $\sik$ if and only if $|i_j - i_k| = 1$.

We are now ready to introduce the second class of patterns that we
wish to avoid.  Say that $w$ is \textit{hexagon-avoiding} if it avoids
each of the patterns in
\begin{equation}
\begin{split}\label{fourhex}
  \{&[4,6,7,1,8,2,3,5],[4,6,7,8,1,2,3,5],\\
    &[5,6,7,1,8,2,3,4],[5,6,7,8,1,2,3,4]\}.
\end{split}
\end{equation}
If we set
\begin{equation}
  \label{hexu}
  u = s_3s_2s_1s_5s_4s_3s_2s_6s_5s_4s_3s_7s_6s_5,
\end{equation}
then the permutations in \eqref{fourhex} correspond to
$u,us_4,s_4u,s_4us_4$.
\begin{figure}[htbp]
  \begin{center}
    \leavevmode
    \input{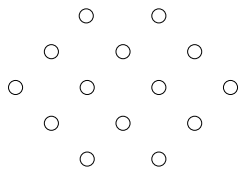}
    \caption{$\heap(u)$ for $u$ as in \eqref{hexu}.}
    \label{fig:hexeg}
  \end{center}
\end{figure}

The heap of any hexagon-avoiding permutation must not contain the
hexagon in \figref{fig:hexeg}.  Permutations that are 321-avoiding and
hexagon-avoiding (321-hexagon-avoiding) can, in fact, be characterized
as those for which no reduced expression contains a substring of
either of the forms
\begin{equation}
\begin{split}
  \label{hexuj}
  u^j = &s_{j+3}s_{j+2}s_{j+1}s_{j+5}s_{j+4}s_{j+3}s_{j+2}\cdot\\
  &s_{j+6}s_{j+5}s_{j+4}s_{j+3}s_{j+7}s_{j+6}s_{j+5}
  \ \ \ \text{ for any } j \geq 0,
\end{split}
\end{equation}
\begin{equation}
  \label{hexuj2}
s_{j}s_{j\pm 1}s_{j}    \text{ for any } j \geq 1.
\end{equation}

It is this characterization of 321-hexagon-avoiding elements that we
will use in the rest of the paper.  

\begin{REM}
  Computationally, it is much more efficient (polynomial time) to
  recognize 321-hexagon-avoiding patterns via pattern avoidance rather
  than by scanning through all reduced expressions for a particular
  subexpression (exponential time).
\end{REM}

The heaps of 321-avoiding elements have a very important property that will
be exploited in the proof of \thmref{mainthm}.  To develop this
property, it will be useful to define the following two subsets of the
unit integer lattice for each $j$, $1\leq j\leq r$:
\begin{align*}
  \text{ The \textit{lower cone}: }\bcone(j) &= 
  \{(i_{j} + \alpha,\rank(\sij) - \beta)\in\bbZ^2: |\alpha| \leq \beta\}.\\    
  \text{ The \textit{upper cone}: }\ucone(j) &= 
  \{(i_{j} + \alpha,\rank(\sij) + \beta)\in\bbZ^2: |\alpha| \leq \beta\}.    
\end{align*}
The \textit{boundary} of $\bcone(j)$ (or $\ucone(j)$) corresponds to the
points in this cone where $|\alpha| = |\beta|$ (see \figref{fig:coneeg}).

\begin{figure}[htbp]
  \begin{center}
    \leavevmode
    \input{coneeg.pstex_t}
    \caption{$\heap(u)$ overlaid with $\bcone(6)$ and $\ucone(6)$.
    The white nodes are in $\heap(u)$.  The black nodes are in one of
    the cones, but not in $\heap(u)$.}
    \label{fig:coneeg}
  \end{center}
\end{figure}

The following lemma yields a very nice property of 321-avoiding
permutations.  In \rmkref{latconrmk}, we interpret this result visually
in terms of $\heap(w)$.

\begin{lem}[Lateral Convexity]\label{latcon}
  Label the generators of $\mathfrak{S}_n$ such that $s_is_j = s_js_i$ if and
  only if $|i-j| > 2$ (the standard labeling).  Then $w\in \mathfrak{S}_n$ is
  321-avoiding if and only if any two occurrences of some
  $s_i$ in a reduced expression for $w$ are separated by both an
  $s_{i-1}$ and an $s_{i+1}$.
\end{lem}

\begin{REM}\label{latconrmk}
  \lemref{latcon} can be rephrased as follows. Suppose that $w =
  s_{i_1}\cdots s_{i_r}$ is 321-avoiding and $\pt(j),\ \pt(k)\in
  \heap(w)$ with $\rank(\sij) < \rank(\sik)$.  Suppose further that
  for each $m\in[i_j,i_k]$ (if $i_j \leq i_k$) or $m\in [i_k,i_j]$ (if
  $i_j > i_k$), there is a point $(m,\rank(s_{i_l}))\in \bcone(\sik)\cap
  \ucone(\sij) \cap \heap(w)$ for some $l,\ j\leq l \leq k$.  Then the
  entire diamond $\bcone(\sik)\cap\ucone(\sij)$ is contained in
  $\heap(w)$.  This is illustrated in \figref{fig:pyrheapeg}.  This
  interpretation relies on Lateral Convexity, that $w$ is
  321-avoiding, and the ``coalescing'' performed in the embedding that
  defines $\heap(w)$.

  \begin{figure}[htbp]
    \begin{center}
      \leavevmode
      \input{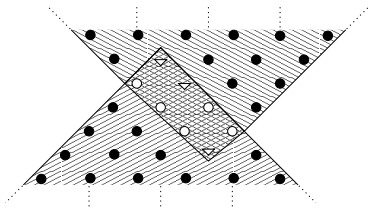}
      \caption{If it is known that the triangular nodes are in $\heap(w)$,
      then \lemref{latcon} tells us that all the white circles are also
      in $\heap(w)$.}
      \label{fig:pyrheapeg}
    \end{center}
  \end{figure}
\end{REM}

\begin{proof}[Proof of Lateral Convexity]
  Suppose $w\in \mathfrak{S}_n$ is 321-avoiding.  Choose a reduced
  expression for $w$ for which a pair of $s_i$'s is as close together
  as possible for some $i$.  These two copies of $s_i$ must be
  separated by at least one of $s_{i\pm 1}$, otherwise our expression
  would not be reduced.  But then our reduced expression looks like
  $u_1 s_i u_2 s_{i\pm 1} u_3 s_i u_4$ where $l(w) = 3 + \sum_{j=1}^4
  l(u_j)$.  If $s_i u_2 = u_2 s_i$ and $u_3 s_i = s_i u_3$, then
  $w$ has a reduced expression $u_1u_2s_is_{i\pm 1}s_iu_3u_4$.  Such a $w$ 
  is not 321-avoiding, which is a contradiction.
  So either $u_2$ or $u_3$ must contain $s_{i\mp 1}$.
  
  For the reverse implication, suppose that every two copies of the
  same generator $s_i$ in some reduced expression for $w$ are
  separated by both an $s_{i-1}$ and an $s_{i+1}$.  It is a theorem of
  Tits \cite{T}, that any two reduced expressions for
  $w\in\mathfrak{S}_n$ can be obtained from each other by a sequence of
  moves of the following two types:
  \begin{align}
    C_1:\ \ &s_is_j = s_js_i, \text{ if } |i-j|>1,\\
    C_2:\ \ &s_is_js_i = s_js_is_j, \text{ if } i = j\pm 1.
  \end{align}
  But, under our hypothesis, we are never able to apply a $C_2$ move
  for such a $w$.  So all reduced expressions for $w$ must be
  obtainable by a sequence of $C_1$ moves.  Hence, $w$ is
  321-avoiding.
\end{proof}

\section{Deodhar's Framework}\label{sec:frame}

For 321-hexagon-avoiding permutations, we will give an explicit
combinatorial formula for the Kazhdan-Lusztig polynomials.  This will
be done in a framework developed by Deodhar \cite{Deodhar90} (using
slightly different notation).  The necessary concepts are reviewed in
this section.

Our construction of the Kazhdan-Lusztig polynomials will be in terms
of subexpressions of a fixed reduced expression $\br = s_{i_1}\cdots
s_{i_r}$.  To this end, we define a \textit{mask} $\bsig$ (associated
to $\br$) to be any binary word $(\sigma_1,\ldots,\sigma_r)$ in the
alphabet $\{0,1\}$.  Set $\bsig[j] \defeq
(\sigma_1,\ldots,\sigma_j)$ for $1\leq j \leq r$.  (So $\bsig =
\bsig[r]$.)  We'll use the notation
\begin{equation}
  \sij^{\sigma_j} = 
  \begin{cases}
    \sij,& \text{ if }\sigma_j = 1,\\
    1,& \text{ if }\sigma_j = 0.
  \end{cases}
\end{equation}
Hence, $w^{\bsig[j]} \defeq s_{i_1}^{\sigma_1} \cdots \sij^{\sigma_j}$
is a (not necessarily reduced) subexpression of $w$.  Let
$\pi(w^{\bsig[j]})$ denote the corresponding element of $\mathfrak{S}_n$.  $\cpw$
will denote the set of ($2^r$ possible) masks of $\br$.  Note that
$\cpw$ can be viewed as the power set of $\{1,\ldots,r\}$.  Finally,
for $x\in \mathfrak{S}_n$, set $\cpxw \subset \cpw$ to be the subset consisting
of those masks $\bsig$ such that $\pi(w^\bsig) = x$.

Define the \textit{defect set} $\cd(\bsig)$ of the fixed reduced expression
$\br$ and associated mask $\bsig$ to be 
\begin{equation}
  \cd(\bsig) = \{j: 2 \leq
  j \leq n,\ l(\pi(w^{\bsig[j-1]})\cdot\sij) < l(\pi(w^{\bsig[j-1]}))\}.  
\end{equation}
Note that $j$'s membership in $\cd(\bsig)$ is independent of
$\sigma_k$ for $k \geq j$.  The elements of $\cd(\bsig)$ are simply
called \textit{defects} (of the mask $\bsig$).

\begin{exmp}
  Let $w = s_3s_2s_1s_4s_3s_2s_5s_4s_3$, $\bsig =
  (1,1,0,1,0,1,0,1,0)$.  Then $w^\bsig = w^{\bsig[9]} =
  s_3s_2s_4s_2s_4$, $\pi(w^\bsig) = s_3$, and $\cd(\bsig) =
  \{6,8,9\}$.  If $x = s_1s_3s_5$, then 
  \begin{equation}
  \begin{split}
    \cpxw = \{\bsig' &= (0,0,1,0,0,0,1,0,1),\\\bsig'' &= (0,0,1,0,1,0,1,0,0),\\
          \bsig''' &= (1,0,1,0,0,0,1,0,0),\\\bsig'''' &= (1,0,1,0,1,0,1,0,1)\}.
  \end{split}
\end{equation}
So, $\cd(\bsig') = \emptyset$, $\cd(\bsig'') = \{9\}$, 
  $\cd(\bsig''') = \{5,9\}$, and $\cd(\bsig'''') = \{5\}$. 
\end{exmp}

Deodhar, in \cite[Lemma 4.1, Definition 4.2, Proposition
4.5]{Deodhar90}, gives a more combinatorial characterization of the
Kazhdan-Lusztig polynomials.  Specifically, he proves that one can
always find a subset $\mathcal{S} \subseteq \cpw$ that yields the
Kazhdan-Lusztig polynomials.  This is an amazing result.  However, in
general, the procedure to find this subset $\mathcal{S}$ is somewhat
complicated.  But we can restrict our attention to the case where
$\mathcal{S} = \cpw$.  In this case, Deodhar's result can be
translated as follows:

\begin{thm}\label{admiss}
  Let $W$ be any finite Weyl group and $\br$ be a reduced expression
  for some $w\in W$.  Set
  \begin{equation}
    \pxw = \sum_{\bsig\in \cpxw} q^{|\cd(\bsig)|}.
  \end{equation}
  If $\deg \pxw \leq \frac{1}{2}\left(l(w)-l(x)-1\right)$ for all
  $x\in W$, then $\pxw$ is the Kazhdan-Lusztig polynomial $P_{x,w}$
  for all $x\in W$.
\end{thm}

Most of the content of \thmref{admiss} is that the $\pxw$ satisfy a
recursive formula equivalent to \thmref{klthm}.

\section{The Defect Graph}\label{sec:graph}

The purpose of the defect graph is to furnish us with a simple
criterion for ensuring that $|\cd(\bsig)| \leq \frac{1}{2}(l(w) -
l(\pi(w^\bsig)) - 1)$ as required by \thmref{admiss}.  However, it is
advantageous to first rephrase this inequality in another language.
So again we introduce some notation.  Partition the defect set
$\cd(\bsig) = \cd^0(\bsig) \cup \cd^1(\bsig)$ where
$\cd^\epsilon(\bsig)$ consists of those $j\in \cd(\bsig)$ for which
$\sigma_j = \epsilon\in \{0,1\}$.  Let $\br[j] \defeq
s_{i_1}\cdots s_{i_j}$ for $1\leq j \leq r$.  Also, set $d_j(\bsig)
\defeq |\cd(\bsig[j])|$, $d(\bsig) \defeq |\cd(\bsig)|$, $x[j] \defeq
\pi(w^{\bsig[j]})$ and $w[j] \defeq \pi(\br[j])$.  Finally, set
\begin{equation}
  \Dbrj \defeq \frac{l(w[j]) - l(x[j]) - 1}{2} - |\cd(\bsig[j])|.  
\end{equation}
We write $\Dbr$ for $\Dbrr$.  Having $\Dbr \geq 0$ implies that the
inequality in \thmref{admiss} holds.  The defect graph will
allow us to show that a condition equivalent to $\Dbr \geq 0$, stated
in the following lemma, holds whenever $w$ is 321-hexagon-avoiding.

\begin{lem}\label{zeroineqlem}
  Let $\br = s_{i_1}\cdots s_{i_r}$ be a reduced expression for some
  $w\in \mathfrak{S}_n$.  Suppose $\bsig=(\sigma_1,\ldots,\sigma_r) \in \cpw$
  with $\pi(w^\bsig) \neq w$.  Then $\Dbr \geq 0$ if and only if
  \begin{equation}\label{zeroineq}
    \left(\text{\# of 0's in } \{\sigma_1,\ldots,\sigma_r\}\right)
    \geq 2\cdot |\cd^0(\bsig)| + 1.
  \end{equation}
\end{lem}
\begin{proof}
  Let $k$ be the smallest index for which $\sigma_k = 0$.  Such a $k$
  must exist by our stipulation that $\pi(w^\bsig) \neq w$.  Consider
  the sequence $w[k],w[k+1],\ldots$.  Since $s_{i_1}\cdots s_{i_k}$ is
  reduced, $\cd(\bsig[k]) = \emptyset$.  Hence, $\Dbrk = 0$.  We now
  investigate the differences $\Dbrj - \Dbrjm$ for $j > k$.  There are
  four possibilities (note that in each case, $l(w[j]) = l(w[j-1]) +
  1$):
  \begin{enumerate}
  \item $j\not\in \cd(\bsig)$, $\sigma_j = 1$.
    Then $d_j(\bsig) = d_{j-1}(\bsig)$, $\ l(x[j]) =
    l(x[j-1]) + 1$.\\  So $\Dbrj - \Dbrjm = 0$.
  \item $j\not\in \cd(\bsig)$, $\sigma_j = 0$.
    Then $d_j(\bsig) = d_{j-1}(\bsig)$, $\ l(x[j]) =
    l(x[j-1])$.\\  So $\Dbrj - \Dbrjm = 1/2$.
  \item $j\in \cd(\bsig)$, $\sigma_j = 1$.
    Then $d_j(\bsig) = d_{j-1}(\bsig) + 1$, $\ l(x[j]) =
    l(x[j-1]) - 1$.\\  So $\Dbrj - \Dbrjm = 0$.
  \item $j\in \cd(\bsig)$, $\sigma_j = 0$.
    Then $d_j(\bsig) = d_{j-1}(\bsig) + 1$, $\ l(x[j]) =
    l(x[j-1])$.\\  So $\Dbrj - \Dbrjm = -1/2$.
  \end{enumerate}
  So, the only cases we need to consider are the second and the
  fourth.  From this it follows that for each $j>k$,
  \begin{equation}
    \Dbrj \geq 0\quad \Longleftrightarrow \quad
    \text{\# of 0's in } \{\sigma_{k+1},\ldots,\sigma_j\} 
    \geq 2 \cdot |\cd^0(\bsig[j])|.
  \end{equation}
  The conclusion of the lemma follows by induction upon setting $j=r$.
\end{proof}

Recall that we need to show that \eqref{zeroineq} is satisfied for
321-hexagon-avoiding permutations for any choice of reduced
expression.  To do this, we define a graph $G_\bsig$ whose vertices
are in one-to-one correspondence with the defects of $\cd^0(\bsig)$.
In \lemsref{rook},\ref{diamond}, and \ref{shapelem} we develop some
technical results relating the shape of $\heap(w)$ to the shape of
$G_\bsig$.  Then in \propref{gtree} we show that $G_\bsig$ is a forest
if $w$ is 321-hexagon-avoiding.  The proof of this Proposition is rather
intricate and is given as a ``proof by picture.''  Finally, in
\secref{sec:pf} we conclude by a simple combinatorial argument that if
$G_\bsig$ is a forest, then \eqref{zeroineq} is satisfied.

\begin{figure}[htbp]
  \begin{center}
    \leavevmode
    \input{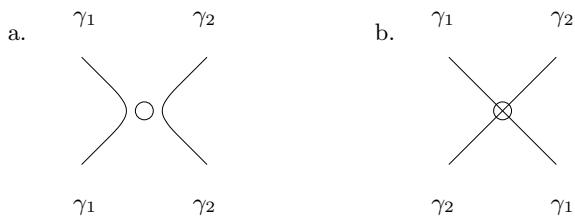}
    \caption{Overlay of string diagram corresponding to some $\bsig$ on $\heap(w)$.}
    \label{fig:bounce}
  \end{center}
\end{figure}

The edges of $G_\bsig$ will depend on how the various defects and
zeros in $\bsig$ are intertwined.  To measure this intertwining, we
overlay strings on $\heap(w)$.  In particular, we will overlay the
lines $y = \pm x + C$ for $C\in \mathbb{Z}$.  At each point $\pt(j)$
of our heap we will move these strings according to the following
rule: If $\sigma_j = 0$, then ``bounce'' the strings as in
\figref{fig:bounce}.a.  If $\sigma_j = 1$, then ``cross'' the strings
as in \figref{fig:bounce}.b.

\begin{figure}[htbp]
  \begin{center}
    \leavevmode
    \input{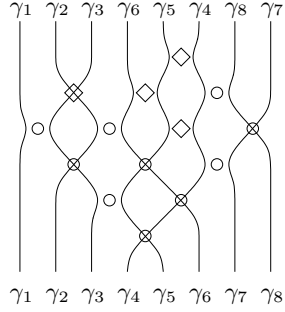}
    \caption{$\heap(w)$ overlaid with a string diagram for the reduced 
      expression $\br = s_4s_3s_2s_1s_5s_4s_3s_2s_6s_5s_4s_7s_6s_5$
      and $\bsig = (1,0,1,0,1,1,0,1,0,0,0,1,0,0)$.  Note that
      $\pi(w^\bsig) = s_4s_5s_4s_7$, giving the permutation
      $[1,2,3,6,5,4,8,7]$.  The defects are represented by diamonds.
      As an illustration of our terminology regarding strings, note
      that $\gamma_4\gamma_7$ meet at $\pt(9)$ (for our reduced
      expression $\br$).  And $\gamma_6$ encounters $\pt(j)$ for $j\in
      \{5,6,7,11\}$ (also for $\br$).}
    \label{fig:stringeg}
  \end{center}
\end{figure}

In either case, $\gamma_1$ and $\gamma_2$ are said to \textit{meet} at
$\pt(j)$ and each of $\gamma_1,\gamma_2$ is said to \textit{encounter}
$\pt(j)$.  If we number the strings from left to right along the
bottom of our heap, reading the order of the strings at the top gives
the permutation $\pi(w^\bsig)$.  \figref{fig:stringeg} gives an example.
\begin{REM}\label{havetocross}
  In the heap model, defects occur when two strings meet that have
  previously crossed an odd number of times.
\end{REM}

\begin{REM}
  In our diagrams, we make the following conventions.  First, every
  diamond point is known to be a defect.  Second, white nodes are known
  to be in our heap.  Third, the inclusion of black nodes within the
  heap is undetermined at the time the picture is first referenced.
\end{REM}

Suppose $j\in \cd(\bsig)$.  For the strings meeting at $\pt(j)$ to
have previously crossed, they both need to have changed direction at
some point (see \figref{fig:critzero}).  Formally, there must be $a,b$
with $1 \leq a \neq b < j$ and $\alpha,\beta > 0$ such that
$(i_j,\rank(\sij)) = (i_a+\alpha,\rank(s_{i_a})+\alpha) =
(i_b-\beta,\rank(s_{i_b})+\beta)$ where $\sigma_{i_a} = \sigma_{i_b} =
0$.  Otherwise, the strings meeting at $(i_j,\rank(\sij))$ could not
have previously crossed.

Choose $a,b$ as above and as large as possible.  Call $\lcz(j) = \pt(a)$
the \textit{left critical zero} and $\rcz(j) = \pt(b)$ the
\textit{right critical zero} of $j$ (or of $\pt(j)$).  In terms of the
heap, the left and right critical zeros ($\lcz(j)$ and $\rcz(j)$) are
the closest zeros to $\pt(j)$ on the boundary of $\bcone(j)$.

Now, for $j\in \cd^0(\bsig)$, $\{\lcz(j), \rcz(j), \pt(j)\}$ are the
\textit{critical zeros} of $j$.  For this reason, we will sometimes
refer to $\pt(j)$ as the \textit{middle critical zero} of $j$ (denoted
$\mcz(j)$).  A point $\pt(j)$ is \textit{shared} if $\pt(j)$ is a
critical zero for two separate defects.

\begin{figure}[htbp]
  \begin{center}
    \leavevmode \input{critzero.pstex_t}
    \caption{Heap showing necessity of existence of 0's on
      boundary of $\bcone(j)$ when $j\in \cd(\bsig)$.}
    \label{fig:critzero}
  \end{center}
\end{figure}

\begin{figure}[htbp]
  \begin{center}
    \leavevmode \input{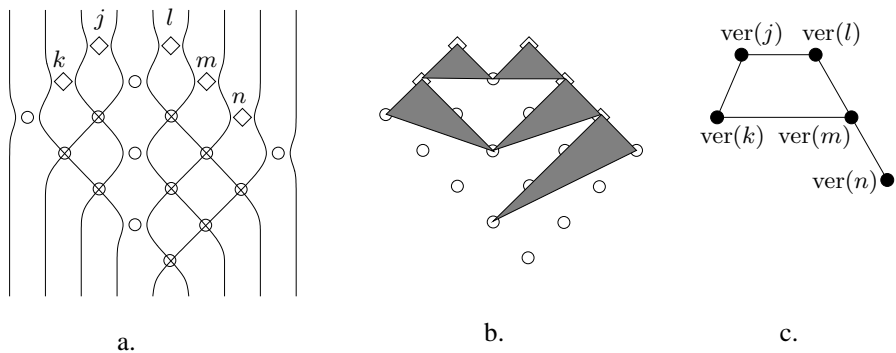}
     \caption{In a., we depict the permutation
         $w = [6,7,8,1,9,2,3,4,5]$ along with the mask $\bsig =
         (1,0,1,1,0,1,1,0,1,0,1,1,1,0,0,0,0,0,0)$.        
       In b., $G_\bsig$ is graphically overlaid on $\heap(w)$.  The
       critical zeros correspond to the corners of the triangles.  In
       c., we have an abstract realization of the graph.}
    \label{fig:grapheg}
  \end{center}
\end{figure}

There is one final construct we will need to prove \thmref{mainthm}.
Define a graph $G_\bsig$ associated to $\bsig$ as follows.  Let the
vertex set of $G_\bsig$ be $\{\ver(j)\}_{j\in \cd^0(\bsig)}$.  The
edge set consists of those $(\ver(j),\ver(k))$ for which
\begin{equation}
  \{\lcz(j),\rcz(j),\mcz(j)\} \cap \{\lcz(k),\rcz(k),\mcz(k)\} \neq \emptyset.
\end{equation}
In \figref{fig:grapheg}, we give an example of a heap along with its
associated graph $G_\bsig$.


The key fact we need in the proof of \eqref{zeroineq} is that
$G_\bsig$ does not contain any cycles.  Before proving this fact in
\propref{gtree}, we first introduce some lemmas that illuminate the
structure of $G_\bsig$.  The first two lemmas are easy and stated only
for reference.  The second and third give criteria for $\heap(w)$ to
contain a hexagon.

\begin{lem}\label{rook}
  Suppose $w$ is 321-avoiding and $k,l\in\cd^0(\bsig)$ with $\pt(l) =
  \lcz(k)$.  Then $\pt(l) + (1,-3)\in\heap(w)$.  Similarly, if $\pt(l)
  = \rcz(k)$, then $\pt(l) - (1,3)\in\heap(w)$.  (See, for example,
  \figref{fig:trip}.c.)
\end{lem}

\begin{lem}\label{diamond}
  Let $w$ be a 321-avoiding permutation and $\pt(h),\pt(k)\in\heap(w)$ with
  $\pt(h) \in \bcone(\pt(k) - (0,6))$.  If $h$ and $k$ are encountered
  by a common string, then $\heap(w)$ contains a hexagon.  (See, for
  example, \figref{fig:trip}.c.)
\end{lem}

\begin{lem}\label{shapelem}
  Let $w\in \mathfrak{S}_n$ be 321-avoiding. $\heap(w)$ contains a
  hexagon if any
  of the following three situations are met:\\
  1.\ The point $\lcz(r) = \pt(m) = \rcz(l)$ with
  $m,r,l\in\cd^0(\bsig)$.  (See \figref{fig:trip}.a.)\\
  2.\ The string $\gamma$ encounters three distinct strings
  $\gamma_1,\gamma_2,\gamma_3$ at defects $l,k,m\in\cd^0(\bsig)$,
  respectively.  Furthermore, $\pt(m) = \rcz(l), \pt(l)=\lcz(k)$ and
  $\pt(m)$ is on the boundary of $\bcone(\pt(k) - (0,2))$.  (See
  \figref{fig:trip}.b.)\\
  3.\ We have $\pt(l) = \lcz(k)$, $\pt(r)=\rcz(k)$ and $k,l,r\in\cd^0(\bsig)$.
  (See \figref{fig:hat}.)
\end{lem}
Parts 1 and 3 of \lemref{shapelem} tell us that any three defects in
a $\vee$-shape or a $\wedge$-shape imply that our heap has a hexagon.
Part 2 of \lemref{shapelem} tells us that, under certain conditions,
if one string encounters three defects, then we also have a hexagon.

\begin{proof}[Proof of part 1.]
  A picture is given in \figref{fig:trip}.a.  The claim follows
  immediately from Lateral Convexity by applying \lemref{rook} to the
  pairs $\pt(l),\pt(m)$ and $\pt(r),\pt(m)$.

  \begin{figure}[htbp]
    \begin{center}
      \leavevmode \input{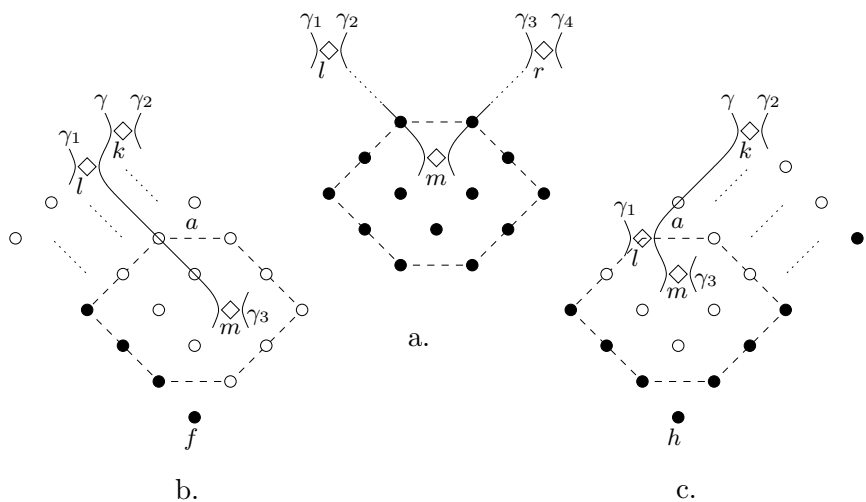}
      \caption{Illustration a. shows the situation of
        \lemref{shapelem}.1.  Illustrations b.,c. refer to
        \lemref{shapelem}.2.  In these latter two pictures, it is
        possible that $\pt(a) = \pt(k)$.}
      \label{fig:trip}
    \end{center}
  \end{figure}

  \noindent
  \emph{Proof of part 2.} First consider the case where $\pt(m) =
  \pt(k) + (\delta,-2-\delta)$ for $\delta \geq 1$.  This is
  illustrated in \figref{fig:trip}.b.  By \lemref{rook}, $\pt(f) =
  \pt(m)-(1,3)$ is in $\heap(w)$.  Since $\pt(f)\in\bcone(\pt(k) +
  (\delta-1,-5-\delta))$, $\heap(w)$ contains the indicated hexagon by
  \lemref{diamond}.
  
  Alternatively, we can have $\pt(m) = \pt(k) - (\delta,2+\delta)$ for
  $\delta \geq 0$.  This is illustrated in \figref{fig:trip}.c.  Recall
  that the $\gamma_i$ are assumed to be distinct.  So, starting at
  $\pt(m) - (1,1)$, $\gamma$ must move down to the right at least
  twice (to cross $\gamma_2$ and $\gamma_3$), and move down to the
  left at least once (to cross $\gamma_1$).  Hence, the lowest of the
  three crossings $\gamma\gamma_i$ must occur in $\bcone(\pt(h)) =
  \bcone(\pt(m) - (0,4)) = \bcone(\pt(a)-(\delta,6+\delta)$.  By
  \lemref{diamond}, $\heap(w)$ must therefore contain a hexagon.

  \begin{figure}[htbp]
    \begin{center}
      \leavevmode \input{hat2.pstex_t}
      \caption{Situation of \lemref{shapelem}.3.}
      \label{fig:hat}
    \end{center}
  \end{figure}

  \noindent
  \emph{Proof of part 3.} 
  By \lemref{rook}, in order to avoid a hexagon in $\heap(w)$, we need
  at least one of $\pt(l),\pt(r)$ to be a distance of exactly
  $\sqrt{2}$ from $\pt(k)$.
  
  Suppose first that both $\pt(l) = \pt(k) - (1,1)$ and $\pt(r) =
  \pt(k) + (1,-1)$.  Then we are in the situation of
  \figref{fig:hat}.a.  Note that if $\sigma_a = 0$ then $a\in
  \cd^0(\bsig)$ and we can appeal to \lemref{shapelem}.1.  So we can
  consider only the case where there is a crossing at $\pt(a)$.  If
  $\gamma$ is either $\gamma_1$ or $\gamma_3$, then it still needs to
  cross a string currently to its right (either $\gamma_2$ or
  $\gamma_4$, respectively).  This can only happen in $\bcone(f)$.
  The only alternative is that $\gamma = \gamma''$.  But then
  $\gamma_1\gamma_2$ cannot cross until $\bcone(f)$.  Either way,
  $\pt(f)\in \heap(w)$.  Arguing analogously with $\gamma'$, we see that
  $\pt(g)\in \heap(w)$.  So $\heap(w)$ contains a hexagon.
    
  Now suppose that only one of $\pt(l),\pt(r)$ is a distance of
  $\sqrt{2}$ away from $\pt(k)$.  Without loss of generality, we
  assume that this point is $\pt(l)$.  We argue depending on whether
  or not $\pt(r)\in\ucone(\pt(m) + (0,2))$ where $\pt(m) = \rcz(l)$.
  
  Assume first that $\pt(r)\in\ucone(\pt(m)+(0,2))$.  We are in the
  situation of \figref{fig:hat}.b.  Since $\pt(r) \neq \pt(k) +
  (1,-1)$, $\pt(m) = \pt(k) + (\delta,-2-\delta)$ for some $\delta
  \geq 2$.  Hence, in order to avoid a hexagon, we must have
  $\gamma_1\gamma_2$ cross as shown.  But then it is easily seen
  that the crossing $\gamma_3\gamma_4$ must occur in $\bcone(h)$.
  This ensures that $\heap(w)$ contains the indicated hexagon.
  
  If $\pt(r)\not\in\ucone(\pt(m)+(0,2))$, then we are in the situation
  of \figref{fig:hat}.c.  Since $\gamma_2$ must go left once below
  $\pt(m) - (1,1)$ (to cross $\gamma_1$) and $\gamma_3$ must go right
  once (to cross $\gamma_4$), we see that the lowest of the crossings
  $\gamma\gamma_i$ must occur in $\bcone(h)$.  If $\pt(r) \neq
  \pt(a)$, then by \lemref{diamond}, $\heap(w)$ contains a hexagon.
  If $\pt(r) = \pt(a)$, then we need the additional fact that $\pt(m)
  \neq \pt(k) - (0,2)$ to ensure that $\pt(h)\in \bcone(\pt(k) - (0,6))$.
  But this follows from the assumption that
  $\pt(r)$ is not at a distance of $\sqrt{2}$ from $\pt(k)$.
\end{proof}

\begin{prop}\label{gtree}
  If $w$ is 321-hexagon-avoiding and $\bsig \in \cpw$, then $G_\bsig$ is
  a forest.
\end{prop}
\begin{proof}
  Assume that $G_\bsig$ is not a forest --- i.e., $G_\bsig$ contains a
  cycle.  We will assume that $w$ is 321-avoiding and show that if
  $G_\bsig$ contains a cycle then $\heap(w)$ contains a hexagon.  Note
  that since $w$ is 321-avoiding, \lemref{latcon} (Lateral Convexity)
  holds.

  Let $V\subset \cd^0(\bsig)$ be a minimal subset such that the
  subgraph $G_\bsig'$ of $G_\bsig$ spanned by $V$ is a cycle.  Hence,
  for each $p\in V$, $\ver(p)\in G_\bsig'$ has degree at least 2.
  Choose $C\in\bbZ$ as large as possible such that $\pt(j)$ is on the
  line $y = x + C$ for some $j\in V$.  Now choose $l\in V$ to be
  minimal among such $j$.  By choice of $V$, $\pt(m) = \rcz(l)$
  must be shared and we must have $\pt(l) = \lcz(k)$ for some
  $k\in V$.  So our heap looks like \figref{fig:findk0}.a.
  \begin{figure}[htbp]
    \begin{center}
      \leavevmode \input{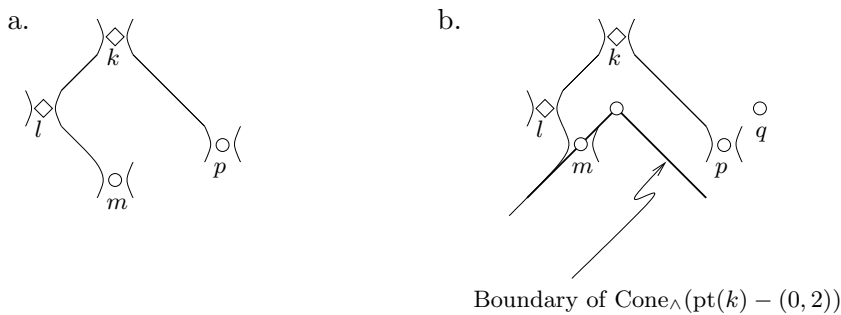}
      \caption{Configuration of $\heap(w)$.  Recall that diamond nodes
        are known defects and white nodes are known to be in
        $\heap(w)$ .}
      \label{fig:findk0}
    \end{center}
  \end{figure}
  
  In the discussion that follows, ``shared'' should be interpreted in
  the context of $G_\bsig'$.
  
  Since $V$ is minimal, either $\pt(k) = \lcz(u)$ for some $u\in V$,
  or $\pt(p) = \rcz(k)$ is shared.  In the first case, $\pt(p) +
  (1,1)$ must be in $\heap(w)$ by Lateral Convexity.  Consider the
  second case --- where $\pt(p)$ is shared.  By \lemref{shapelem}.3,
  $p\not\in V$.  So $\pt(p) = \lcz(r)$ for some $r\in V$.  So in both
  cases, we have the following fact which we state for reference.

  \begin{fact}\label{qfact}
    If $\pt(p) = \rcz(k)$, then $\pt(q) = \pt(p) + (1,1)\in \heap(w)$.
  \end{fact}

  Two other simple facts we state for reference are the following.

  \begin{fact}\label{fact:0} By Lateral Convexity, any point
    encountered by a string that still needs to cross below that point
    must be in the heap (after pushing together connected components).
    For example, if $j\in\cd(\bsig)$, then $pt(j) - (0,2)$ must be in
    the heap.
  \end{fact}
  \begin{fact}\label{mfact}
    Recall that $\pt(m)$ is defined as right critical zero of the left
    critical zero of $\pt(k)$ (see \figref{fig:findk0}.b).  If
    $\heap(w)$ does not contain a hexagon, then the point $m$ must lie
    along the boundary of $\bcone(\pt(k)-(0,2))$.
  \end{fact}

  We now show that, regardless of the characteristics of $m$ (i.e.,
  values of $i_m$, $\rank(m)$, and whether or not $m\in\cd(\bsig)$),
  $\heap(w)$ must contain a hexagon.  Suppose that $m\in V$.  By
  \lemref{shapelem}.2, the only way this can happen is if the other
  string encountering $\pt(m)$ is $\gamma_3$.  Since $V$ is minimal,
  we then need either $\lcz(m)$ or $\rcz(m)$ shared.  Consider
  \figref{mpos}.  Suppose $\pt(n) = \lcz(m)$ is shared.  By choice of
  $\pt(k)$ on the line $y = x + C$, this implies that $n\in V$.  But then
  by \lemref{rook}, $\pt(h)\in\heap(w)$.  Then by Lateral Convexity,
  $\pt(e)\in\heap(w)$.  The alternative is that $\rcz(m)$ is shared.
  Again, this implies that $\pt(e)\in\heap(w)$.  Since $\pt(q)\in\heap(w)$
  by \factref{qfact}, $\heap(w)$ contains a hexagon.

  \begin{figure}[htbp]
    \begin{center}
      \leavevmode \input{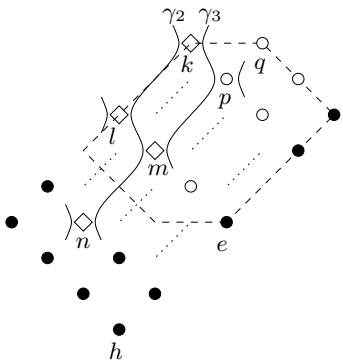}
     \caption{This figure depicts the case where $\pt(m)$ is not the
       left critical zero of another defect in $V$.}
      \label{mpos}
    \end{center}
  \end{figure}
  So we can assume that $m\not\in V$.  But by choice of $l$, $\pt(m)$
  must be shared.  This implies that $\pt(m) = \lcz(r)$ for some $r\in
  V$.  We now argue that $\heap(w)$ must contain a hexagon according
  to the position of $\pt(m)$ relative to $\pt(k)$.
  
  \textbf{Case I:} $\pt(m) = \pt(k) - (\delta,2+\delta)$ for $\delta
  \geq 0$.  There are three cases to consider.  \figref{upright}.a.
  depicts the first.  Here, $\gamma$ and $\gamma_3$ both encounter
  $\pt(r)$.  Since $V$ is minimal, either $\rcz(r)$ or $\pt(r)$ must
  be shared.  By choice of our line $y = x + C$ and the fact that
  $p\not\in V$, we see that, in fact, $\rcz(r)$ must be shared.  But
  then $\pt(b)\in\heap(w)$.  Since $\pt(q)\in\heap(w)$, $\heap(w)$
  contains the indicated hexagon.
  
  \begin{figure}[htbp]
    \begin{center}
      \leavevmode \input{upright.pstex_t}
      \caption{}
      \label{upright}
    \end{center}
  \end{figure}

  The second alternative is that $\pt(r) \in \bcone(\pt(k))$ but
  $\gamma_3$ does not encounter $\gamma$ along any of the nodes
  between $\pt(m)$.  This is depicted in \figref{upright}.b.  If
  $\sigma_c = 0$, then $\gamma_2\gamma_3$ must cross in $\bcone(g)$.
  If $\sigma_c = 1$ , then $\gamma\gamma'$ must cross in $\bcone(e)$.
  In either case, $\heap(w)$ must contain the indicated hexagon.
  
  The third possibility is that $\pt(r)\not\in\bcone(\pt(k))$
  (\figref{upright}.c).  In fact, this is the only possibility for
  $\pt(r)$ when $\delta = 0$.  Here we see that the path of $\gamma_3$
  must be as shown in order to avoid $\bcone(g)$.  But then
  $\gamma_4\gamma_5$ cannot cross until $\bcone(e)$.  So we have the
  indicated hexagon.

  \textbf{Case II:} $\pt(m) = \pt(k) + (\delta,-2-\delta)$ for some
  $\delta \geq 1$.  The situation is depicted in \figref{caseIII}.a.
  For both $\gamma_1\gamma_2$ and $\gamma_2\gamma_3$ to cross outside
  of $\bcone(h)$, we need $\gamma_2\gamma_3$ to cross in $\ucone(m)$.
  This is shown in \figref{caseIII}.b.  We mention three additional
  assertions we have made in \figref{caseIII}.b.  First, $\gamma_1$
  must cross $\gamma_2$ as shown in \figref{caseIII}.b in order to
  avoid having $\heap(w)$ contain a hexagon.  Second, $\pt(q)\in \heap(w)$
  by \factref{qfact}.  Third, since $\rcz(m)$ must be shared,
  $\pt(e)\in\heap(w)$ as shown.  So, by Lateral Convexity, $\heap(w)$
  contains the hexagon indicated in \figref{caseIII}.b.  (It is
  possible that $\pt(a) = \pt(p)$ or $\pt(a) = \pt(k)$, but this does
  not change our conclusion.)
  \begin{figure}[htbp] 
    \begin{center} \leavevmode
      \input{casethree.pstex_t} 
      \caption{} 
      \label{caseIII} 
    \end{center}
  \end{figure}
\end{proof}

\section{Proof of \thmref{mainthm}}\label{sec:pf}

We present one remaining needed lemma and then the proof of
\thmref{mainthm}.  

In the following lemma, we let $\br = s_{i_1} \cdots s_{i_r}$ be a
reduced expression for $w$ and set $s = s_{i_r}$.  Then let $\br/s$
denote the truncated reduced expression $s_{i_1}\cdots s_{i_{r-1}}$
for $ws$.

\begin{lem}\label{rec}
  Let $s\in S$, $ws < w$.  Then
  \begin{equation}
    \pxw = q^{c_s(x)} \pxws + q^{1-c_s(x)}\pxsws,\label{pxrec}
  \end{equation}
\end{lem}
\begin{equation}
  \text{ where } c_s(x) =
  \begin{cases}
    1, & \text{ if } xs < x,\\
    0, & \text{ if } xs > x.\\
   \end{cases}
\end{equation}

\begin{proof}
  Partition $\cpxw = \cpzero \dot{\cup} \cpone$ where $\cpeps$
  consists of all masks in $\cpxw$ ending in $\epsilon$ for
  $\epsilon\in \{0,1\}$.  There are natural bijections $\cpone \approx
  \cpxsws$ and $\cpzero \approx \cpxws$ given by $\bsig \mapsto
  \bsig[r-1]$.  So, to prove the lemma, we need only compare
  $|\cd(\bsig)|$ to $|\cd(\bsig[r-1])|$.
  
  If $\bsig\in \cpzero$, then $\bsig[r-1]\in \cpxws$.  In this case,
  if $xs > x$ ($c_s(x) = 0$), then $r\not\in\cd(\bsig)$, so
  $|\cd(\bsig[r-1])| = |\cd(\bsig)|$.  Alternatively, if $xs < x$
  ($c_s(x) = 1$), then $\cd(\bsig) = \cd(\bsig[r-1]) \cup \{r\}$ and
  $|\cd(\bsig)| = |\cd(\bsig[r-1])| + 1$.  This accounts for the first
  term in \eqref{pxrec}.
  
  Since $c_s(xs) = 1 - c_s(x)$, proof of the second term in
  \eqref{pxrec} reduces to the above case.
\end{proof}

\begin{proof}[Proof of \thmref{mainthm}]
  \textit{\ref{itone} $\Longrightarrow$ \ref{ittwo}:}
  
  Assume $w$ is 321-hexagon-avoiding.  We need to show that the
  $\pxw$ are the Kazhdan-Lusztig polynomials.  
  
  Now, every $j\in \cd^0(\bsig)$ has three critical zeros.
  Furthermore, by \lemref{shapelem}, no point is a critical zero for 3
  distinct defects.  So the number of edges in $G_\bsig$ equals the
  number of shared critical zeros.  Hence,
  \begin{align}
    \text{\# of 0's in } \{\sigma_1,\ldots,\sigma_r\}
    &\geq \text{\# of critical zeros in }
    \{\sigma_1,\ldots,\sigma_r\}\label{otone}\\
    &= 3\cdot |\cd^0(\bsig)| - \text{ (\# of edges in } G_\bsig).\label{ottwo}
  \end{align}
  Now, by \propref{gtree}, $G_\bsig$ is a forest with $|\cd^0(\bsig)|$
  vertices.  Hence, $G_\bsig$ has at most $|\cd^0(\bsig)| - 1$ edges
  (see, e.g., \cite{bollobas}).  Hence,
  \begin{equation}
    \text{\# of 0's in } 
    \{\sigma_1,\ldots,\sigma_r\}\geq 2\cdot |\cd^0(\bsig)| + 1.
  \end{equation}
  So by \lemref{zeroineqlem}, $\Dbr \geq 0$.  Therefore the inequality
$|\cd(\bsig)| \leq \frac{1}{2}(l(w) - l(\pi(w^\bsig))-1)$ holds.  Now
apply \thmref{admiss}, from which it follows that $\pxw = P_{x,w}$ for
all $x\in W$.

\textit{\ref{ittwo} $\Longrightarrow$ \ref{itone}:}\\
We shall prove (not \ref{itone}) $\Longrightarrow$ (not \ref{ittwo}).
Assume $w$ is not 321-avoiding.  We can find a reduced expression for
$w$ of the form $\br = v s_i s_{i\pm1} s_i v'$ with $l(w) = l(v) +
l(v') + 3$.  Set
  \begin{equation}
    \bsig = (\overset{l(v)}{\overbrace{1,\ldots,1}},1,0,0,
             \overset{l(v')}{\overbrace{1,\ldots,1}}).
  \end{equation}
  Then $|\cd^0(\bsig)| = 1$ and $|\{j: \sigma_j = 0\}| = 2$.  By
  \lemref{zeroineqlem}, $\Dbr < 0$.  So $\pxw$ does not satisfy the
  properties of the $P_{x,w}$ listed in \thmref{klthm}.
  
  Now assume $w$ is 321-avoiding but not hexagon avoiding.  Then we
  can write $w = vu^jv'$ where $u^j$ as in \secref{sec:prelim} and
  $l(w) = l(v) + l(v') + 14$.  Set
  \begin{equation}\label{bigmask}
    \bsig = (\overset{l(v)}{\overbrace{1,\ldots,1}},
    1,1,0,1,0,1,0,1,1,0,0,0,0,0,
    \overset{l(v')}{\overbrace{1,\ldots,1}}).
  \end{equation}
  The mask $\bsig$ is depicted graphically in \figref{fig:stairs}.
  Then $|\cd^0(\bsig)| = 4$ and $|\{j: \sigma_j = 0\}| = 8$.  By
  \lemref{zeroineqlem}, $\Dbr < 0$.  So $\pxw$ does not satisfy the
  properties of the $P_{x,w}$ listed in \thmref{klthm}.

  \begin{figure}[htbp]
    \begin{center}
      \leavevmode
      \input{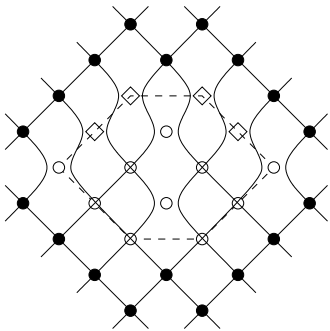}
      \caption{Heap view of mask in \eqref{bigmask}.  The black nodes
      are not known to be in the heap.}
      \label{fig:stairs}
    \end{center}
  \end{figure}

  \textit{\ref{ittwo} $\Longrightarrow$ \ref{itthree}:}\\
  We first appeal to a result of Kazhdan and Lusztig relating the
  intersection Poincar\'e polynomial of the Schubert variety $X_w$ to
  the Kazhdan-Lusztig polynomials $P_{x,w}$ (\cite[Corollary~4.9]{K-L2}):
  \begin{equation}  \label{icpoly}
    \sum_i \dim(\IH^{2i}(X_w)) q^i  = \sum_{x \leq w} q^{l(x)} P_{x,w}(q).
  \end{equation}
  Now, we are assuming that $\pxw = P_{x,w}$ for all $x\in \mathfrak{S}_n$.
  So we need only show that 
  \begin{equation}
    \sum_{x \leq w} q^{l(x)} \pxw = (1+q)^{l(w)}.
  \end{equation}
  We proceed by induction, the result being obvious for $l(w) = 1$.
  Choose an $s\in S$ such that $ws < w$.  From \cite[Lemma 7.4]{Hum},
  we know that: 
  \begin{equation}
    \label{humref}
    \text{ If } ws < w, \text{ then } x \leq w \Longleftrightarrow  xs \leq w.
  \end{equation}
  Using \eqref{humref}, along with \lemref{rec}, we can write
  \begin{align}
    \sum_{x\leq w} q^{l(x)} \pxw &= 
    \sum_{x\leq w,\ x < xs} q^{l(x)} \pxw + q^{l(xs)} \pxsw\label{parallel}\\ 
    &= (1+q)\sum_{x\leq w,\ x < xs} q^{l(x)}
      \left(\pxws + q \pxsws\right)\\
    &= (1+q)\sumsb{x\leq w,\ x < xs} q^{l(x)}\pxws + q^{l(xs)} \pxsws.
  \end{align}
  If $\pxws \neq 0$, then $x\leq ws$, so this becomes
  \begin{align}
    &= (1+q)\sum_{x\leq ws} q^{l(x)} \pxws\\
    &= (1+q)(1+q)^{l(ws)} = (1+q)^{l(w)}.
  \end{align}
  The last line is by the induction hypothesis.

  \textit{\ref{itthree} $\Longrightarrow$ \ref{ittwo}:}\\
  Deodhar \cite{Deodhar90} proves that for any Weyl group $W$, we can
  always find a subset $\mathcal{S}\subset \cpw$ such that
  \begin{equation}
    \sum_{\substack{\bsig\in \mathcal{S}\\\pi(w^\bsig) = x}}
    q^{|\cd(\bsig)|} = P_{x,w}
  \end{equation}
  for all $x,w\in W$.  (More generally, he shows that such an
  $\mathcal{S}$ exists when the coefficients of $P_{x,w}$ are already
  known to be non-negative.  Due to their interpretation in terms of
  dimensions of intersection cohomology groups, this is known for any
  Weyl group.)
  
  Hence, for such an $\mathcal{S}$, we have the following string of
  equalities:
  \begin{equation}
    \begin{split}
      (1+q)^{l(w)} = \sum_i \dim(\IH^{2i}(X_w))q^i &= \sum_{x\leq w}
      q^{l(x)}P_{x,w}\\ &= \sum_{\substack{\bsig\in
          \mathcal{S}\\\pi(w^\bsig) = x}} q^{l(x)}q^{|\cd(\bsig)|}.
  \end{split}
\end{equation}
  Setting $q = 1$, we find that $2^{l(w)} = |\mathcal{S}|$.  But then
  $S = \cpw$.  So $\pxw = P_{x,w}$ for all $x,w\in \mathfrak{S}_n$.

  \textit{\ref{ittwo} $\Longleftrightarrow$ \ref{itfour}:}\\
  This follows from Deodhar
  \cite[Proposition 3.5 and Corollary 4.8]{Deodhar90}.

  \textit{\ref{itthree} $\Longrightarrow$ \ref{itfive}:}\\
  This is the content of Deodhar \cite[Proposition 3.9]{Deodhar90}.  

  \textit{\ref{itfive} $\Longrightarrow$ \ref{itsix}:}\\
  This is a standard result on small resolutions.  See, for example,
  \cite[Section 6.5]{kirwan}.
  
  \textit{\ref{itsix} $\Longrightarrow$ \ref{itthree}:}\\
  Recall that $Y$ denotes the Bott-Samelson resolution of $X_w$
  (corresponding to some reduced expression $\br$).  By
  \cite[Proposition 4.2]{B-S},
  \begin{equation}    
    \sum_i \dim(H^{2i}(Y))q^i = (1+q)^{l(w)}.\label{bseq}
  \end{equation}
  We are assuming that $H_*(Y) \cong \IH_*(X_w)$.  By Poincar\'e
  duality, we know that $H^{2i}(Y) \cong \IH^{2i}(X_w)$.
  Combining \eqref{bseq} with this isomorphism yields
  \begin{equation}    
    \sum_i \dim(\IH^{2i}(X_w))q^i = (1+q)^{l(w)}
  \end{equation}\label{bseqtwo}
  as desired.

  This completes the proof of the \thmref{mainthm}.
\end{proof}

\begin{cor}\label{pxone}
  If $w = s_{i_1}\cdots s_{i_r}$ with $i_1,\cdots,i_r$ all distinct,
  then $P_{x,w} = 1$ for all $x\leq w$.
\end{cor}

\section{A Conjecture of Haiman and a Generalization}\label{sec:conj}

Define $q$-Fibonacci numbers by $F_n(q) = F_{n-1}(q) + qF_{n-2}(q)$
where $F_n(q) = 0$ if $n < 0$ and $F_0(q) = F_1(q) = 1$.
\thmref{mainthm} gives us a simple proof of the following conjecture
of Haiman (\cite[Conjecture 7.18]{Brenti98}):
\begin{cor}\label{cor:haiman}
  Let $w_{k,l} \in \mathfrak{S}_n$ have reduced expression 
  \begin{equation}
    \br = s_ks_{k-1}s_{k+1}s_k\cdots s_ls_{l-1} \in \mathfrak{S}_n,\  
    2 \leq k < l < n.    
  \end{equation}
Then $P_{e,w_{k,l}} = F_{l-k+1}(q)$.
\end{cor}
Recently, Brenti-Simion \cite{Bren-Sim} have independently proved this
conjecture and generalized it to a class of elements that are not
321-hexagon-avoiding.  In fact, the corollary can be generalized to
apply to any 321-hexagon-avoiding element for which no generator
appears more than twice.

\begin{figure}[htbp]
  \begin{center}
    \leavevmode \input{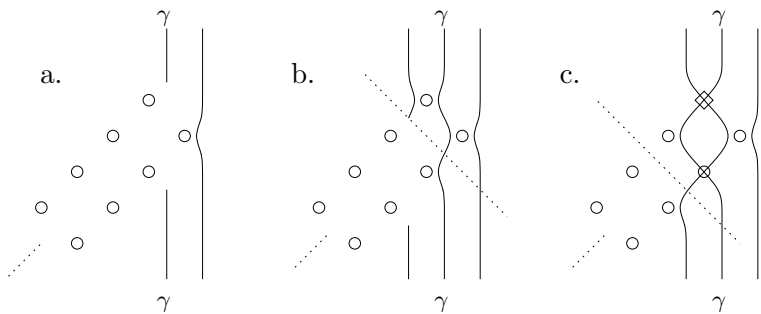}
    \caption{It is clear that $\gamma$ must remain in its column in
      order for $\pi(w^\bsig) = e$.  This is shown in a.  Diagrams b.
      and c. show the only two possibilities for the path of
      $\gamma$.}
    \label{haimanfig}
  \end{center}
\end{figure}

\begin{proof}
  As a permutation, 
  \begin{equation}
    w = [1,2,\cdots,k-2,k+1,\cdots,l+1,k-1,k,l+2,\cdots,n].    
  \end{equation}
  This is easily seen to be 321-hexagon-avoiding.  So by
  \thmref{mainthm}, $P_e(\br) = P_{e,w_{k,l}}$.
  
  The claim is true for $l=k$.  The proof is by induction.  The
  situation of the general case is illustrated in \figref{haimanfig}
  for some $\bsig\in\cpew$.  Let $r = l(w)$.  In
  \figref{haimanfig}.b, no new defect is introduced by $\gamma$, so
  $|\cd(\bsig)| = |\cd(\bsig[r-2])|$.  In \figref{haimanfig}.c, we
  have $|\cd(\bsig)| = |\cd(\bsig[r-4])| + 1$.  The claim follows by
  the induction hypothesis.
\end{proof}

We give below the generalization where
$\heap(w)$ is a $3 \times (l-k+1)$ diamond rather than a $2 \times
(l-k+1)$ diamond.

\begin{thm}\label{3xnthm}
  Suppose $v_{k,l}\in \mathfrak{S}_n$ has reduced expression
  \begin{equation}
    \label{3nwd}
    \br = s_l s_{l+1} s_{l+2} s_{l-1} s_l s_{l+1} \cdots
    s_k s_{k+1} s_{k+2}
  \end{equation}
  for some $1\leq k \leq l < n-2$.  Then $P_{e,v_{k,l}} \in \bbZ[q]$ is
  given by the coefficient of $z^{l-k+1}$ in the generating function
  \begin{equation*}\label{3xnrec}
    G_e(z) = \frac{-1 + q^2z^2 + q^3z^3}
    {(1 + qz + q^2z^2)(-1 + z + qz + qz^2 + q^2z^2 + q^2z^3 - q^4z^4)}.
  \end{equation*}
\end{thm}

\begin{figure}[htbp]
  \begin{center}
    \leavevmode \input{3xn.pstex_t}
    \caption{}
    \label{3xnfig}
  \end{center}
\end{figure}

\begin{proof}
  We only sketch the proof.  We see that $v_{k,l}$ is clearly
  321-hexagon-avoiding, so by \thmref{mainthm}, $P_{x,w} = P_{\br}$.
  The idea is to use recursion on $n = l-k$.  From \figref{3xnfig}, it
  is easy to see that $P_{e,v_{k,l}} = P_{e,v_{k+1,l}} + q
  P_{s_{k+1},v_{k+1,l}} + q P_{s_{k+2},v_{k+1,l}} + q^2
  P_{s_{k+1}s_{k+2},v_{k+1,l}}$.  Similar recurrences can be found for
  $P_{x,v_{k,l}}$ where $x\in \mathfrak{S}_4$.  Solving these
  recurrences for $P_{e,v_{k,l}}$ yields \eqref{3xnrec}.
\end{proof}

\section{Singular Loci of 321-hexagon-avoiding Elements}\label{sec:loci}

The Schubert variety $X_w$ is said to be \textit{singular} at a point
$x\leq w$ (or, more properly, on the Schubert cell $C_x \subset X_w$)
if the Zariski tangent space to $X_w$ at $x$ has dimension strictly
greater than $l(x)$.  The set of singular points forms a lower order
ideal in the Bruhat-Chevalley order (\cite{BLak}).  We define $\xsing$
to consist of the maximal elements (under this Bruhat-Chevalley order)
of the set $\{x\in\mathfrak{S}_n : x\leq w \text{ and } x \text{
  singular}\}$.

The following theorem gives a complete description of $\xsing$ when
$w$ is 321-hexagon-avoiding.  In fact, this proves a conjecture of
Lakshmibai and Sandhya \cite{Lak-San} in this special case.

\begin{thm}\label{maxsing}
  Let $w \in \mathfrak{S}_n$ be 321-hexagon-avoiding (hence $\heap(w)$
  is well-defined).  Then every diamond with vertices $(x,y),
  (x-\alpha,y-\alpha),(x+\beta,y-\beta),(x-\alpha+\beta,y-\alpha-\beta),
  \alpha,\beta > 0$ in the heap determines an element in $\xsing$.
  More explicitly, let
  \begin{multline}
    \label{locitrip}
    T = \{(j,k,l)\in \bbZ^3: 1\leq j,k,l\leq r,\ \pt(j) = \pt(k) -
    (\alpha,\alpha),\\\pt(l) = \pt(k) + (\beta,-\beta) \text{ for
    some } \alpha,\beta > 0,\\ \text{ and }
    \bcone(j) \cap \bcone(l)\cap \heap(w) \neq \emptyset\}
  \end{multline}
  and
  \begin{equation}
  \begin{split}
    \label{}
    \Sigma = \{\bsig \in \cpw: (j,k,l)\in T,\ 
    &\sigma_j = \sigma_k = \sigma_l = 0,\\
    &\text{ and } \sigma_m = 1 \text{ for } m\neq j,k,l\}.
  \end{split}
\end{equation}
Then the maximal singular locus $\xsing$ of $X_w$ is given by
  $\xsing = \{\pi(\bsig): \bsig\in \Sigma\}$.
\end{thm}

\begin{proof}
  It has been proved by Deodhar \cite{Deodhar85} that for
  $W=\mathfrak{S}_n$ and $v\leq w$, $X_w$ is smooth on the Schubert
  cell $C_v$ if and only if $P_{v,w} = 1$.  By \thmref{mainthm},
  $\cpxw = P_{x,w}$ for every $x\in \mathfrak{S}_n$.  So to show that
  $X_w$ is singular, we need only show that $\cpw$ contains a mask of
  positive defect.
  
  Let $\bsig\in\Sigma$ correspond to $(j,k,l)\in T$.  Since every
  defect must have two critical zeros (in addition to the defect
  itself), $l(w) - l(\pi(\bsig)) = 3$.  Lateral Convexity tells us
  that if $l(w) - l(\pi(\bsig)) < 3$ for some other $\bsig\in\cpw$,
  then $|\cd(\bsig)| = 0$.  So for $\bsig\in\Sigma$, if $X_w$ is
  singular at $C_{\pi(\bsig)}$, $\pi(\bsig)$ is maximally singular.
  Now, the conditions in \eqref{locitrip} imply that $k\in\cd(\bsig)$.
  By \thmref{mainthm}, this implies that $P_{\pi(\bsig),w} \neq 1$.
  So $\{\pi(\bsig): \bsig\in \Sigma\} \subseteq \xsing$.
  
  The only fact that remains to be checked is that if $y$ is a
  singular point of $X_w$, then $y\leq \pi(\bsig)$ for some $\bsig\in
  \Sigma$.  So pick some $\bsig\in\cpyw$ with $|\cd(\bsig)| \geq 1$.
  Choose $b\in\cd(\bsig)$ and suppose $\pt(a) = \lcz(b)$ and $\pt(c) =
  \rcz(b)$.  Now define a mask $\bsig'$ by setting
  \begin{equation}
    \sigma_m' = 
    \begin{cases}
      1, & m \not\in \{a,b,c\},\\
      0, & m \in \{a,b,c\}.
    \end{cases}
  \end{equation}
  Using the characterization of Bruhat-Chevalley order in terms of
  subexpressions (see, e.g., \cite{Hum}), it is easily checked that
  $\pi(\bsig) \leq \pi(\bsig')$.  Since $\bsig'$ is in $\Sigma$, we
  are done.
\end{proof}

\begin{cor}
  For $w$ 321-hexagon-avoiding, each element of $\xsing$ has
  codimension 3 in $X_w$.
\end{cor}

\begin{exmp}
  Here we give an example of calculating the singular locus as in
  \thmref{maxsing}.  We have set $w = s_2s_1s_5s_4s_3s_2s_6s_5s_4s_3$.
  \figref{fig:locuseg} illustrates the eight different points in the
  maximal singular locus of $X_w$.  Namely, \begin{align*} \xsing =
  \{&[3,1,6,2,7,4,5],[1,6,3,2,7,4,5],[3,1,6,4,2,7,5],\\
  &[3,1,6,5,2,4,7],[1,3,7,2,6,4,5],[3,2,6,1,4,7,5],\\
  &[3,2,6,1,5,4,7],[3,4,6,1,2,5,7]\}.  \end{align*}
  \begin{figure}[htbp] \begin{center} \leavevmode
  \input{locuseg.pstex_t} \caption{} \label{fig:locuseg} \end{center}
  \end{figure}
\end{exmp}

\begin{exmp}
  For $v_{1,4}$ as in \thmref{3xnthm}, $|\xsing| = 18$.
\end{exmp}

\begin{REM}
  Let $w = [w(1),\dots,w(n)]$.  A result of Lakshmibai and Sandhya
  \cite[Theorem 1]{Lak-San} is that $X_w$ is nonsingular if an only if
  $w$ avoids [3,4,1,2] and [4,2,3,1]. It is shown in \cite{carrell94}
  that $X_w$ is non-singular precisely when $P_{e,w} = 1$.  So from
  \thmref{mainthm} and \corref{pxone}, we see that if $w$ is
  321-hexagon-avoiding and $X_w$ is singular, then we must be able to
  find a [3,4,1,2]-sequence in $w$.
\end{REM}

\section{Example and Enumeration of 321-hexagon-avoiding Elements}\label{sec:enum}

The following table lists both the number of 321-avoiding elements in
$\mathfrak{S}_n$ and the number of 321-hexagon-avoiding elements in
$\mathfrak{S}_n$ for $7\leq n\leq 13$ (these numbers are equal for
$n\leq 6$).  The number of 321-hexagon-avoiding elements has been
calculated by computer.  The number of 321-avoiding elements is
well-known to be given by the Catalan numbers (see, e.g.,
\cite{BJS,Knuth,SS}).

\begin{table}[h]
   \begin{tabular}[c]{rrrrrrrrr}\hline
      n & 7 & 8 & 9 & 10 & 11 & 12 & 13\\\hline\hline
      321-avoiding         & 429 & 1430 & 4862 & 16796 & 
      58786 & 208012 & 742900\\\hline 
      321-hexagon-avoiding & 429 & 1426 & 4806 & 16329 & 
      55740 & 190787 & 654044\\\hline 
\end{tabular}
\caption{Number of 321-hexagon-avoiding elements in $\mathfrak{S}_n$.}
\end{table}

Below we give an example showing the use of \thmref{mainthm} for
calculating $P_{x,w}$.
\begin{exmp}
  Here we calculate $P_{x,w}$ for $w = s_2s_1s_3s_2s_4s_3$.  As a
  permutation, $w$ is $[3,4,5,1,2]$, which is clearly
  321-hexagon-avoiding.    (Note that $w = w_{2,4}$ in the sense of
  \corref{cor:haiman}.)  It is a result of Deodhar that for each $x
  \leq w$, there exists a unique mask in $\cpxw$ of defect 0.  The
  following table lists all of the $\bsig \in \cpw$ for which
  $|\cd(\bsig)| > 0$.  For this $w$, all of these masks happen to have
  $|\cd(\bsig)| = 1$.
  \begin{table}[h]
  \begin{tabular}{llll}\hline
    $s_2s_1s_3s_2s_4s_3$ & $\pi(w^\bsig)$  &  
    $s_2s_1s_3s_2s_4s_3$ & $\pi(w^\bsig)$ \\\hline
     0\;\:0\;\:1\;\:0\;\:0\;\:1\;\: &  $e    $
     &   1\;\:1\;\:1\;\:0\;\:0\;\:1\;\: & $ s_2s_1  $\\
     1\;\:0\;\:0\;\:1\;\:0\;\:0\;\: &  $e    $
     &   1\;\:0\;\:0\;\:0\;\:0\;\:1\;\: & $ s_2s_3   $\\
     0\;\:1\;\:1\;\:0\;\:0\;\:1\;\: &  $s_1 $
     &   1\;\:0\;\:1\;\:0\;\:0\;\:0\;\: & $ s_2s_3   $\\
     1\;\:0\;\:1\;\:0\;\:0\;\:1\;\: &  $s_2  $
     &   1\;\:0\;\:1\;\:1\;\:0\;\:1\;\: & $ s_3s_2   $\\
     1\;\:0\;\:0\;\:0\;\:0\;\:0\;\: &  $s_2  $
     &   1\;\:0\;\:0\;\:1\;\:1\;\:1\;\: & $ s_4s_3   $\\
     1\;\:0\;\:0\;\:1\;\:0\;\:1\;\: &  $s_3  $
     &   1\;\:0\;\:0\;\:0\;\:1\;\:0\;\: & $ s_2s_4   $\\
     0\;\:0\;\:1\;\:0\;\:0\;\:0\;\: &  $s_3  $
     &   1\;\:0\;\:0\;\:0\;\:1\;\:1\;\: & $ s_2s_4s_3$\\
     1\;\:0\;\:0\;\:1\;\:1\;\:0\;\: &  $s_4  $
     &   1\;\:0\;\:1\;\:1\;\:0\;\:0\;\: & $ s_2s_3s_2$\\
     0\;\:1\;\:1\;\:0\;\:0\;\:0\;\: &  $s_1s_3$
     &   1\;\:1\;\:1\;\:0\;\:0\;\:0\;\: & $s_2s_1s_3 $\\\hline
  \end{tabular}
  \caption{Computing $P_{x,w}$ using the defect statistic.}
\end{table}
  Hence, we see that for $x\leq w$,
  \begin{equation}
    P_{x,w} = 
    \begin{cases}
      1+2q, & \text{ if } x \in \{e,s_2,s_3,s_2s_3\},\\
      1+q, & \text{ if } x \in \{s_1,s_4,s_1s_3,s_2s_1,s_3s_2,s_4s_3,\\
          & \qquad \qquad s_2s_4, s_2s_4s_3, s_2s_3s_2, s_2s_1s_3\},\\
      1, & \text{ otherwise. }
    \end{cases}
  \end{equation}
\end{exmp}

\section{Acknowledgements}
We would like to thank Mark Haiman, Allen Knutson, and Rodica Simion.
The second author would like to personally thank the computer program
\textsf{GAP/CHEVIE}.

\bibliography{hexagon}
\end{document}